\documentclass[letterpaper,10pt]{article}

\usepackage{graphicx}
\usepackage{amssymb}
\usepackage{amsmath}
\usepackage{latexsym}
\usepackage[english]{babel}

\newtheorem{teo}{Theorem}
\newtheorem{lem}{Lemma}
\newtheorem{rem}{Remark}
\newtheorem{coro}{Corollary}
\newtheorem{propo}{Proposition}

\title{A global Carleman estimate in a transmission wave equation and application to a one-measurement inverse problem.} 

\date{}

\author{Lucie Baudouin\footnote{e-mail: {\tt baudouin@laas.fr}}\\
{\it\footnotesize LAAS - CNRS; Universit\'e de Toulouse; 7, avenue du Colonel Roche, F-31077 Toulouse, France.}\\ 
Alberto Mercado and Axel Osses\footnote{e-mail: {\tt amercado@dim.uchile.cl, axosses@dim.uchile.cl}}\\
{\it\footnotesize  DIM-CMM, Universidad de Chile, Casilla 170/3 - Correo 3, Santiago, Chile.}}

\begin{document}

\maketitle

\begin{abstract}
We consider a transmission wave equation in two embedded domains in ${\mathbb R}^2$, 
where the speed is $a_1>0$ in the inner domain and $a_2>0$ in the outer domain. 
We prove a global Carleman inequality for this problem under the hypothesis
that the inner domain is strictly convex and $a_1>a_2$. 
As a consequence of this inequality, uniqueness and Lipschitz stability 
are obtained for the inverse problem of retrieving a stationary potential for the 
wave equation with Dirichlet data and discontinuous principal coefficient from a 
single time-dependent Neumann boundary measurement. 
\end{abstract}



\def\RR{\mathbb{R}}
\def\ZZ{\mathbb{Z}}
\def\NN{\mathbb{N}}
\def\1{1 \!\! I}
\def\cD{{\mathcal{D}}}
\def\cU{{\mathcal{U}}}
\def\cC{{\mathcal{C}}}
\def\O{{\mathcal{O}}}
\newcommand{\dn}[1]{\frac{\partial #1}{\partial \nu}}
\newcommand{\dns}[2]{\frac{\partial #1}{\partial \nu_{#2}}}
\newcommand\E{\varepsilon}
\newcommand{\Bo}{B_{\varepsilon}(x_0)}
\newcommand{\Bi}{B_{\varepsilon}(x_1)}
\newcommand{\Bd}{B_{\varepsilon}(x_2)}
\newcommand{\Rr}{\mathbb{R}}
 \newcommand{\om}{\omega}
 \newcommand{\Om}{\Omega}
 \newcommand{\Omo}{\Omega_{x_0}}
  \newcommand{\Gm}{\Gamma}
   \newcommand{\Sg}{\Sigma}
 \newcommand{\la}{\lambda}
 \newcommand{\gm}{\gamma}
 \newcommand{\inS}{\iint_{\Sigma}}
 \newcommand{\Tom}{\omega \times (0,T)}
 \newcommand{\TOm}{\Omega \times (0,T)}
\newcommand{\LOm}{L^2(T_0,T;L^2(\Om))} 
 \newcommand{\LO}{L^2(T_0,T;L^2(\O))}
 \newcommand{\lom}{L^2(\om)}
 \newcommand{\Hom}{H^1_0(\om)}
\newcommand{\InTO}{\int_0^T \int_{\Omega}}
\newcommand{\InTo}{\int_0^T \int_{\omega}}
\newcommand{\fun}[3]{#1 : #2 \longrightarrow #3}
\newcommand{\dem}{\noindent{\em Proof:\ \/}}
\newcommand{\QED}{\hfill{\rule{2mm}{2mm}}}
\font\dsrom=dsrom10 scaled 1000
\def \1{\textrm{\dsrom{1}}}
\newcommand{\tn}{\textnormal}
\newcommand{\ds}{\displaystyle}


\section{Introduction and main results}\label{Intro}
\subsection{Presentation of the problem}
The inverse problem of recovering coefficients from a wave equation with discontinuous coefficients 
from boundary measurements arises naturally in geophysics and more precisely, in seismic prospection
of Earth inner layers~\cite{VanDer}. 

Here we are interested in the case where only one particular measurement is available. This could be important, for instance, in seismic prospection, where data of a single wave that propagates through the Earth is considered.

Consider two embedded domains, where the speed 
coefficients are $a_1>0$ in the inner domain and $a_2>0$ in the outer domain.
Stability of the inverse problem we study here 
is obtained by deriving a global Carleman estimate for the wave equation with 
discontinuous coefficients. We prove this Carleman inequality in the case
the inner domain is strictly convex and the speed 
is monotonically increasing from the outer to the inner layers, i.e. $a_1>a_2$. 
This last situation is, incidentally, the general case into the Earth.

\begin{figure}[ht!]
\begin{center}
\includegraphics[width=10cm]{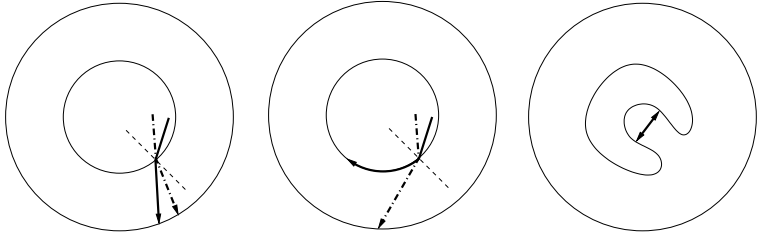}
\caption{\small Two domains with speed coefficients $a_1$ (inner) and $a_2$ (outer). In the first figure
(left),  if $a_1>a_2$ by Snell's law all the inner rays reach the exterior boundary independently of their incident angles. Conversely, in the second figure (center) if $a_1<a_2$ some rays with large incident 
angles remain trapped near the inner interface. The last figure (right) shows a trapped ray into a captive domain.}
\end{center}
\label{Fig1}
\end{figure}

Figure~\ref{Fig1} illustrates the role of these hypothesis and gives some intuition with the 
help of Snell's law. In the case $a_1>a_2$  (see Figure~\ref{Fig1}, left) 
the incident rays coming from the inner domain toward the outer domain become closer to the normal 
at the interface since $\sin(\theta_1)>\sin(\theta_2)$, where $\theta_i$, $i=1,2$ are the corresponding
 incident angles. Therefore, all the rays coming from the inner ball with any incident angle $\theta_1$ 
 in $(-\pi/2,\pi/2)$ succeed in crossing the interface. In the opposite case, when $a_1<a_2$ 
 (see Figure~\ref{Fig1}, center) we have 
$\sin(\theta_1)<\sin(\theta_2)$ and there is a critical incident angle $\theta^*<\pi/2$ such that 
the rays with incident angles $\theta_1$ out of the range $(-\theta^*,\theta^*)$ remain supported 
near the interface and do not reach the outer domain, so this information does not arrive 
at the exterior boundary. Finally, strict convexity of the inner domain avoids trapped rays 
(see Figure~\ref{Fig1}, right).

Global Carleman estimates and the method of Bukhgeim-Klibanov \cite{Bu-Kli}, \cite{B} are especially useful for solving the one measurement inverse problems. It is possible to obtain local Lipschitz stability around the single known solution, provided that this solution is regular enough
and contains enough information \cite{Kli-Mal} (see also \cite{Kl} and \cite{Y}).
Many other related inverse results for hyperbolic equations use the same strategy. A complete list is too long to be given here. To cite some of them see \cite{P-Y} and \cite{Y} where Dirichlet boundary data and Neumann measurements are considered and \cite{ImYam2}, \cite{I-Y} where Neumann boundary data and Dirichlet measurements are studied. These references are all based upon the use of local or global Carleman estimates. Related to this, there are also general pointwise Carleman estimates that
are also useful in similar inverse problems \cite{Fu,XuZhang,Kli-Ti}. 

Recently, global Carleman estimates and applications to one-measurement inverse problems were obtained in the case of variable but still regular coefficients, see \cite{ImYam1} for the isotropic case, and \cite{LasTriYao} and \cite{Bella} for the anisotropic case. It is interesting to note that these authors require a bound on the gradient of the coefficients, so that the idea of approximating discontinuous coefficients by smooth ones is not useful.

\begin{figure}[ht!]
\begin{center}
\includegraphics[height=3.5cm]{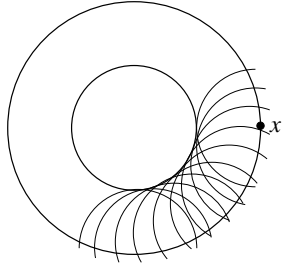}
\caption{\small Recovering the interface as the envelope of circumferences. Each circumference, centered at some point $x$ on the exterior boundary, represents the possible locations of the nearest point of the inner interface where the reflection took place for a given traveltime measured at $x$.}
\end{center}
\label{Fig2}
\end{figure}

There are a number of important works \cite{Rak1,Han,Rak2,Symes,Beylkin}
concerning the same inverse problem in the case that several boundary measurements 
are available. In these cases, it is possible to retrieve speed coefficients and 
even discontinuity interfaces without any restrictive hypothesis of strict convexity or 
speed monotonicity. For instance, one can retrieve the interface by observing the
traveltime reflection of several waves. Indeed, it is well known that the interface can be recovered as the envelope of certain curves as shown in Figure~\ref{Fig2} (see also \cite{VanDer} and the references therein). This method works independently of the sign of $a_1-a_2$ and this explains in part why there are no geometrical or speed monotonic hypotheses for these kind of inverse results. 

Let us now give some insight into the relationship between this work and exact controllability or 
energy decay for the wave equation with discontinuous coefficients. 

First of all, the global Carleman estimate we obtain immediately implies a particular case of a 
well known result of exact controllability for the transmission wave equation \cite{L}. Roughly speaking, 
the result of \cite{L} states that we can control internal waves from the exterior boundary in a layered speed media if the speed is monotonically increasing from  the outer to the inner layers and the inner domain is star shaped, a weaker assumption than strict convexity. Moreover, if the speed monotonicity is inverted, there are non controllable solutions with concentrated energy near
the interface \cite{Ma-Zu}, \cite{Ca-Zu}. 

Secondly, there exist several results about the growth
of the resolvent for the spectral stationary transmission problem, from where it is possible to 
derive the speed of local energy decay for the evolution wave equation with transmission conditions
\cite{Burq}. In the case $a_1>a_2$ and if the inner domain is strictly convex, it has been shown using micro-local analysis \cite{PopVo1,PopVo2} that the speed of the energy decay is exponential if the dimension of space is odd and polynomial otherwise. In the general case, including the cases 
when $a_1<a_2$ or the inner domain is not strictly convex, it has been proved using micro-local analysis and global Carleman estimates for the spectral problem \cite{Bella1} that  the energy 
decays as the inverse of the logarithm of the time.

Notice that we shall only consider here the case of a discontinuous coefficient which 
is constant on each subdomain (i.e. $a_1$ and $a_2$ constants). We will indeed concentrate our discussion on the main difficulty, namely the discontinuity at the interface. However, we could also consider variable coefficients $a_1(x)$ and $a_2(x)$  such that their traces at the interface are constant,   under additional assumptions of boundedness of $\nabla a_j$ similar to those appearing in \cite{ImYam1} (see Remark \ref{Coefvar}).

Finnaly, we note that a global Carleman estimate \cite{Doubo2} has also been obtained for the heat equation with discontinuous coefficients. That work was initially motivated by the study of the exact null controllability of the semilinear wave equation, but the estimate has been recently used to prove
local Lipschitz stability for a one measurement inverse problem for the heat equation with discontinuous coefficients \cite{Bella2}, \cite{Bena}. 

Having introduced the problem, let us now present our main results.

\subsection{Inverse problem}
Let  $\Om$ and  $\Om_1$ be two open subsets of $\RR^2$ with smooth boundaries $\Gm$ and $\Gm_1$. Suppose that $\Om_1$ is simply connected, $\overline{\Om}_1 \subset \Om$ and set
$\Om_2 = \Om \setminus \overline{\Om}_1$. Thus, we have $\partial \Om_2 = \Gm \cup \Gm_1$.
We also set:
\begin{equation} \nonumber
a(x) = \left\{ \begin{array}{cc}
a_1    &x\in \Om_1 \\
a_2    &x\in \Om_2
\end{array}
\right. 
\end{equation}
with $a_j > 0$ for $j = 1,2$.  We consider the following wave equation:
\begin{equation} \label{waveqn}
 \left\{ \begin{array}{rclc}
         u_{tt}  - \textnormal{div}( a(x) \nabla u)  + p(x)u & = & 0 & (x,t)\in\Omega \times (0,T)\\
                              u & = & 0 & (x,t)\in\Gamma \times (0,T)\\
                           u(0) & = & u_0  & x\in\Omega\\
                           u_t(0) & = & u_1 & x\in\Om.
\end{array}
\right. 
\end{equation}

We know that \cite{L-M,Da-Li}
for each $p \in L^{\infty}(\Om)$, $u_0 \in H_0^1(\Om)$ and $u_1 \in L^2(\Om)$, there exists a unique weak solution $u(p)$ of equation (\ref{waveqn}) such that
$$
u(p) \in C([0,T];H_0^1(\Om)), \,u_t\in C([0,T];L^2(\Om)).
$$
We shall prove the well-posedness of the inverse problem consisting of retrieving the potential $p$ involved in equation (\ref{waveqn}), by knowing the flux (the normal derivative) of the solution $u(p)$ of (\ref{waveqn}) on the boundary. 
We will prove uniqueness and stability of the non linear inverse problem characterized by the  non linear application 
\begin{equation}\label{nonlinear}
\left.p\right|_{\Omega}\longmapsto\left.a_2\frac{\partial u}{\partial \nu}\right|_{\Gamma\times(0,T)}.
\end{equation}
More precisely, we will answer the following questions.

\noindent\textbf{Uniqueness} : \\
Does the equality $\frac{\partial u(q)}{\partial\nu}= \frac{\partial u(p)}{\partial \nu}$ on $\Gamma\times(0,T)$ imply $q=p$ on $\Omega$ ?\\

\noindent\textbf{Stability} : \\
Is it possible to estimate $\left.\left(q-p\right)\right|_{\Omega}$ by  $\left.\left(\frac{\partial u(q)}{\partial \nu}- \frac{\partial u(p)}{\partial \nu}\right)\right|_{\Gamma\times(0,T)} $ in suitable norms~?\\

The idea is to reduce the {\it nonlinear inverse problem} 
to some {\it perturbed inverse problem} which will be solved with the help of a 
global Carleman estimate. More precisely, we will give a local answer about the determination of
$p$, working first on the perturbed version of the problem, as shown is Section~3. Assuming that $p\in L^{\infty}$ is a given function, we are concerned with the
stability around $p$. That is to say, $p$ and $u(p)$ are known while $q$ and $u(q)$
are unknown. 

We are able to prove the following result, which states the stability of the inverse problem. 

\begin{teo} \label{ip}
Assume $\Om$ is bounded, $\Om_1\subset\subset\Omega$ is a strictly convex domain 
with boundary $\Gamma_1$ of class $C^3$ and  $a_1>a_2 > 0$. 
There exists $T_0>0$ such that, given $T>T_0$, if  
$p\in L^\infty(\Om)$, $u_0\in H^1_0(\Om)$,  $u_1\in L^2(\Omega)$ and $r  > 0$ satisfy
\begin{itemize}
\item $|u_0(x)| \geq r > 0$ a. e. in $\Om$, and
\item $u(p) \in H^1(0,T;L^{\infty}(\Om))$  
\end{itemize}
then, given a bounded set   $\mathcal U \subset L^\infty(\Om)$,  there exists a constant \\
$C = C(a_1,a_2, \Omega_1,\Omega_2, T, \|p\|_{L^{\infty}(\Om)}, \|u(p)\|_{H^1(L^\infty)}, \mathcal U,r)>0$
such that:
$$\|p - q\|_{L^2(\Om)}  \leq C \left\|a_2\dn{u(p)}  - a_2\dn{u(q)}\right\|_{H^1(0,T;L^2(\Gm))}$$
for all $q \in \mathcal U$, where $u(p)$ and $u(q)$ are the solutions of (\ref{waveqn}) with potential $p$ and $q$, respectively.
\end{teo}

\begin{rem}
In section \ref{PbI} is given an estimate for  $T_0$ in function of $a_1$, $a_2$, $\Om_1$ and $\Om_2$. 
See  Theorem \ref{teolp}.
\end{rem}

Let us remark that as a direct consequence of the \textit{local stability} of Theorem~\ref{ip}
we have the following \textit{global uniqueness} for our inverse problem: 
\begin{coro}
If $u(p)$ and 
$u(q)$ are two solutions of (\ref{waveqn}) for potentials $p$ and $q$
in $L^\infty(\Omega)$ with $u(p)$, $u_0$, $u_1$, $a_1$, $a_2$, $\Omega_1$, $\Omega_2$ and $T$ 
satisfying the hypothesis of Theorem~\ref{ip} and such that 
$\frac{\partial u(q)}{\partial\nu}=\frac{\partial u(p)}{\partial \nu}$ on 
$\Gamma\times(0,T)$ then $p=q$.
\end{coro}

The proof of this result consists of two parts: a global Carleman estimate and the resolution of the
inverse problem and lipschitz stability following the methods introduced in \cite{Bu-Kli} and \cite{Kli-Mal}
which we have already mentioned in the introduction.

\subsection{Carleman estimate}

We introduce here our main result concerning a global Carleman estimate for the solutions 
of problem (\ref{waveqn}) extended to the time interval $(-T,T)$. We set 
$Q=\Omega\times(-T,T)$, $\Sigma= \Gamma\times(-T,T)$, $\Sigma_1= \Gamma_1\times(-T,T)$, 
$Q_j=\Omega_j\times(-T,T)$, $u_j = u \1_{Q_j}$ and
$\nu_j$ the outward unit normal to $\Om_j$, for $j = 1,2$.

We will work with an equivalent formulation of (\ref{waveqn}).
Notice that for each $f \in L^2(Q)$,  $u$ solves the equation
\begin{equation}
u_{tt}  - div(a \nabla u) + pu =  f \quad \mbox{  in  }  Q \hspace{0.5cm}
\end{equation}
if and only if, for each $j \in \{ 1, 2\}$,   $u_j$ solves (see \cite{L})
\begin{equation}
u_{j,tt}  - a_j \Delta u_j + pu_j =  f\1_{Q_j}   \quad  \mbox{  in  }  Q_j \hspace{0.5cm}
\end{equation}
together with the {\it transmission conditions}
\begin{equation}\label{Tr}             
  \left\{ \begin{array}{cl}
         u_1 = u_2   &\textnormal{on }\Sigma_1\\
                    a_1 \dns{u_1}{1} + a_2 \dns{u_2}{2} = 0  &\textnormal{on }\Sigma_1.
          \end{array}
\right. 
\end{equation}

In order to construct a convenient weight function, take $x_0 \in \Om_1$ and for each $x\in\Om\setminus \{x_0\}$ define $\ell(x_0,x) = \{ x_0 + \la(x - x_0) \, : \,  \la \geq 0 \}$.  
Since $\Om_1$ is convex there is exactly one point $y(x)$ such that
\begin{equation}\label{defy}
y(x) \in \Gm_1 \cap \ell(x_0,x).
\end{equation}
We define the function 
$\fun{\rho}{\Om\setminus \{x_0\}}{\RR}^+$ by:
\begin{equation} \label{rho}
\rho(x) = |x_0 - y(x)|. 
\end{equation}
Let $\E > 0$ be such that $\overline B_\E\subset\Omega_1$ (and small enough in a sense
that we will precise later) and let $0<\E_1<\E_2<\E$. Then we consider a cut-off function 
$\eta \in C^{\infty}(\RR)$ such that 
\begin{equation}
0 \leq \eta \leq 1,\qquad
\eta = 0 \textrm{ in } B_{\E_1}(x_0),\qquad
\eta = 1 \textrm{ in } \Omega\setminus \overline B_{\E_2}(x_0).
\end{equation}
For each $j \in \{ 1,2\}$ we take $k$ such that $\{j,k \} = \{1,2\}$ and we define the following
functions in the whole domain $\Omega\times\RR$ 
\begin{equation}\label{phi-j}
\phi_j(x,t)  =  \eta(x) \frac{a_k}{\rho(x)^2}|x - x_0|^2 - \beta t^2 + M_j  
\qquad  (x,t) \in \Om \times\RR,
\end{equation}
where $\beta$, $M_1$ and $M_2$ are positive numbers that will be chosen later. 
Then, the weight function we will use in this work is
\begin{equation}\label{phi}
\hspace{1cm} 
\phi(x,t) = 
\left\{ \begin{array}{lrr}
 \phi_1(x,t)     
 & (x,t) \in \Om_1 \times \RR & \vspace{0,3cm}\\ 
\phi_2(x,t)    
& (x,t) \in \Om_2 \times \RR. &
\end{array}
\right.
\end{equation}
Notice that (see $(c)$ and $(d)$ in Proposition \ref{propphi} below) 
$\phi_1$ and $\phi_2 $ satisfy (\ref{Tr}) if and only if
\begin{equation} \label{H1}  
 M_1 - M_2 = a_1 - a_2. 
\end{equation}
We denote 
$$L = \partial_t^2 - a \Delta~~\textnormal{and}~~E(z) = |z_t|^2 - a |\nabla z|^2.$$ 
As usual, we do the change of variables
\begin{equation}\label{varphi}
\varphi  =  e^{\lambda \phi}  \mbox{ , }  \lambda > 0,  \qquad
w  =  e^{s \varphi} u \mbox{ , }  s > 0,\qquad
P(w)  =   e^{s \varphi} L( e^{- s \varphi} w)
\end{equation}
and after algebraical computations, we split $P(w)$ into three terms as follows:
$$
P(w) = P_1(w) + P_2(w) + R(w),
$$
where for some fixed real number $\gm \in (0,1)$
\begin{eqnarray*}
P_1(w) & = & w_{tt} - a\Delta w + s^2\lambda^2\varphi^2E(\phi)w, \\
P_2(w) & = & (\gm-1)s \lambda \varphi L(\phi) w - s \lambda^2 \varphi E(\phi) w 
- 2s \lambda \varphi (\phi_t w_t - a \nabla \phi \cdot \nabla w), \\
R(w)  & =  & - \gm s \lambda \varphi L(\phi) w.
\end{eqnarray*}
We will write $P^{\phi}$, $P_1^{\phi}$,  $P_2^{\phi}$, etc. if we want to make the 
dependence on $\phi$ explicit. Also,
given $U \subset \RR^2$,  we define the norm in $ H^1(U \times (-T,T))$, given by
\begin{equation} \label{normphi}
\| g \|_{_{U, \varphi}}^2 =
  s \la \int^T_{-T}\int_U (|g_t|^2 +  |\nabla g|^2 )\varphi
 + s^3\la^3 \int_{-T}^T\int_U|g|^2 \varphi^3,
 \end{equation}
and
 \begin{equation} \label{bd_0}
 \Sigma^\phi_+ = \{ (x,t)\in\Sigma : \nabla \phi(x,t) \cdot \nu(x) > 0 \}.
 \end{equation}
Finally, we define the space
\begin{eqnarray*}
X = \{ u \in  L^2(-T,T;L^2(\Om)) & | & Lu_j \in L^2(-T,T;L^2(\Om_j)), \, j=1,2; \, u|_{\Sg} = 0, \\
	&      &    u(\pm T) = u_t(\pm T) = 0, \mbox{ and } u \mbox{ satisfies  (\ref{Tr})} \}.
\end{eqnarray*}

The main global Carleman estimate is the following
\begin{teo} \label{Carleman}
Assume $\Om_1$ is a strictly convex domain of class $C^3$, and $a_1> a_2>0$. 
Let $x_k\in \Om_1$, $k=1,2$ and let
$\phi^k$, $\varphi^k$, $w^k$
be the corresponding functions defined for $x_k$ as we did before for $x_0$
in $(\ref{defy})$, $(\ref{rho})$, $(\ref{phi})$ and $(\ref{varphi})$.  
Let $\nu$ be the unit outward normal to $\Om$.
Then there exists $C>0$, $s_0>0$ and $\la_0 >0$ such that
 \begin{gather}
 \sum_{k=1}^{2} \left(  
 \left\|P_1^{\phi^k}(w^k) \right \|^2_{L^2(Q)} \right.+ \left. \left \|P_2^{\phi^k}(w^k) \right \|^2_{L^2(Q)}
+   \left\|w^k \right\|_{_{Q,\varphi^k}}^2   \right) \nonumber  \\
 \leq C \sum_{k=1}^2 \left(
  \left \|P^{\phi^k}(w^k) \right\|^2_{L^2(Q)}   +  s \la \displaystyle \iint_{\Sigma_+^{\phi^k}} \varphi^k \left|a_2 \dn{w^k} \right|^2  \right)  \label{ineqCarlem}
 \end{gather}
for all $u \in X$,  $\la \geq \la_0$ and $s \geq s_0$.
\end{teo}
Notice that in the right-hand side of $(\ref{ineqCarlem})$ we have the term
$$
\left\|P^{\phi^k}(w^k)\right\|_{L^2(Q)}^2  =  \iint_{Q} e^{2s\varphi^k}\left|Lu\right|^2.$$ 
Since we consider an equation given by the operator $L_p = \partial_{tt} - a \Delta + p$, it is important to note that the same estimate with  right-hand side equal to  $\ds\iint_{Q} e^{2s\varphi^k}\left|L_pu\right|^2$ is also true for all potentials $p$ such that $|p|_{L^{\infty}(\Om)} \leq m$, with $m$ already fixed. Indeed, 
$$
|Lu|^2 \leq 2|L_pu|^2 + 2m^2|u|^2
$$
and taking $s$ large enough,  the left hand side of the Carleman estimate of Theorem~\ref{Carleman} can absorb the term $2Cm^2\ds\iint_{Q} e^{2s\varphi}|u|^2$. That is, we have the following result.
\begin{coro}\label{carlem}
Under the hypothesis and notations of Theorem $\ref{Carleman}$, given $m \in\RR$,  
there exists $C>0$ (depending on $m$), $s_0>0$ and $\la_0 >0$ such that for all $p\in L^{\infty}(\Om)$ with $\|p\|_{L^{\infty}(\Om)} \leq m$ we have 
 \begin{gather}
  \sum_{k=1}^{2} \left(  
 \left\|P_1^{\phi^k}(w^k) \right \|^2_{L^2(Q)}  \right. +  \left.   \left \|P_2^{\phi^k}(w^k) \right \|^2_{L^2(Q)}
+   \left\|w^k \right\|_{_{Q,\varphi^k}}^2   \right)    \nonumber \\
 \leq C \sum_{k=1}^2 \left(
  \iint_Q e^{2s \varphi_k}|L_p u|^2   
   +  s \la  \iint_{\Sigma_+^{\phi^k}} \varphi^k \left|a_2 \dn{w^k} \right|^2  \right)   \label{ineqCarlem_p}
\end{gather}
for all $u \in X$,  $\la \geq \la_0$ and $s \geq s_0$.
\end{coro}

The paper is organized as follows. In Section~\ref{Carlem} we give the proof of Theorem~\ref{Carleman}. 
In Section~\ref{PbI} we apply inequality (\ref{ineqCarlem_p}) to derive the stability of the inverse problem presented in Theorem~\ref{ip}.
The proof of Theorem~\ref{Carleman} is organized in several subsections. 
In Proposition \ref{propphi} we  prove the properties of  $\phi$ which will allow us to use it as a weight in a Carleman estimate. In subsection \ref{list}, with $\phi$ as the weight function,  we develop the $L^2$-product $\left< P_1(w),P_2(w) \right>$ for functions with non-zero boundary values. We prove  inequality (\ref{ineqCarlem}) in subsection (\ref{proofCarlem}).

\section{Proof of the Carleman inequality}\label{Carlem}

\subsection{Weight function}

Here we prove that the function $\phi$ satisfies enough properties for being a 
weight function in a Carleman estimate.

We will use the following notation:
\begin{quote}
$M = M_1 \1_{Q_1} + M_2 \1_{Q_2}$, ~~ $\bar{a} = a_2 \1_{Q_1} + a_1 \1_{Q_2}$, ~~  $c(x) = \ds{\frac{\bar{a}}{\rho(x)^2}}$, \\
$\Om_0 = \Om_1 \cup \Om_2$, ~~ $\Omo = \Om_0 \setminus\overline\Bo$, \\
$Q_0 =\Omega_0 \times(-T,T)$, ~~ $Q_{x_0} = \Omega_{x_0} \times(-T,T)$.
\end{quote}

\begin{propo} \label{propphi}
  If $\Om_1$ is a {\bf strictly convex} domain of class $C^3$,  we can take $\E, \delta > 0$ such that:
\begin{enumerate}
\item[(a)] $|\nabla \phi| \geq \delta > 0$ \hspace{0.5cm} in $Q_{x_0} =  (\Om_1 \cup \Om_2)\setminus \Bo \times (-T,T)$
\item[(b)] $\nabla\phi_1(x,t) \cdot \nu_1(x) \geq \delta > 0$,  \hspace{0.5cm} $ \forall \, (x,t) \in \Sg_1$
\end{enumerate}
where $\nu_1$ is the unit outward normal vector to $\Omega_1$.\\
If additionally  $(\ref{H1})$ is satisfied, we also have: 
\begin{enumerate}
\item[(c)] $\phi_1(x,t) = \phi_2(x,t) = a_2-\beta t^2 + M_1$ \hspace{0.5cm} $\forall \, (x,t) \in \Sg_1$
\item[(d)] $\displaystyle{ a_1  \frac{\partial^{\alpha_1+\alpha_2}}{\partial x_1^{\alpha_1}\partial x_2^{\alpha_2}}  \phi_1(x,t) = 
a_2  \frac{\partial^{\alpha_1+\alpha_2}}{\partial x_1^{\alpha_1}\partial x_2^{\alpha_2}}  \phi_2(x,t)}$ \hspace{0.5cm} for all $ \, (x,t) \in \Sg_1 \, $ and  
$\alpha_1, \alpha_2 \in \NN \cup \{ 0 \}$ with $\alpha_1 + \alpha_2 \leq 3$.
\item[(e)] $\Delta \phi(x,t) \geq 2 c(x)$ \hspace{0.5cm}   $\forall (x,t) \in Q_{x_0} $.
\item[(f)] $D^2(\phi)(X,X) \geq \delta_1 |X|^2$ \hspace{0.5cm}  in $Q_{x_0}   \hspace{0.3cm} \forall X \in \RR^2$ \hspace{0.4cm} for some $\delta_1 > 0$.
\end{enumerate}
\end{propo}

\dem
We have $$\nabla \phi = 2c(x)(x-x_0) + |x-x_0|^2\nabla c(x).$$ 

By definition, $\rho(x)$ (thus also $c(x)$) is constant in the direction of $x-x_0$. Therefore
$$ (x-x_0) \cdot \nabla c(x) = 0 $$
and 
\begin{eqnarray*}
|\nabla \phi|^2 & = & 4c^2(x)|x-x_0|^2 + |x-x_0|^4|\nabla c(x)|^2 \\
		&  \geq &  4c^2(x)|x-x_0|^2 \\
		&  \geq  & 4 \left( \frac{\bar{a}}{\textrm{diam}(\Omega)^2} \right)^2 \E^2  \hspace{0.5cm} \mbox{ in } \Omo
\end{eqnarray*}
and $(a)$ is proved.

Now, it is clear that $\phi_1(x,t) = a_2-\beta t^2 + M_1$ for each $(x,t) \in \Sigma_1$,  so $\Gamma_1 \times \{t\}$ is a level curve of $\phi_1(\cdot, t)$ for each $t \in [-T,T]$. Since 
$\phi_1(x,t) < a_2 - \beta t^2 + M_1 < \phi_1(y,t)$ for any $x \in \Om_1$ and $y \in \Om_2$, we have $\nabla \phi_1 = |\nabla \phi_1| \nu_1 $ on $\Sigma_1$, and thus $(a)$ implies $(b)$.

By definition  $\rho (x) = | x-x_0|$ for all $x \in \Gamma_1$, hence  $(c)$ is simply deduced from $(\ref{H1})$.

Without lost of generality, we can take  $x_0 = 0$. Writing $\rho$ in polar coordinates, $\Gamma_1$ can be parameterized by
\begin{equation} \label{paramrho}
\gamma(\theta) = (\rho(\theta)\cos\theta, \rho(\theta)\sin\theta).
\end{equation}
and  then $\rho$  is a $C^3$ function. If $\ds{D = \frac{\partial^{\alpha_1+\alpha_2}}{\partial x_1^{\alpha_1}\partial x_2^{\alpha_2}} }$ with $\alpha_1 + \alpha_2 \leq 3$, we get
$$\quad a_1 D   \phi_1(x,t) = a_1 a_2 D \left( \frac{|x-x_0|^2}{\rho(x)^2}\right) = a_2 D  \phi_2(x,t)  $$
for all $(x,t) \in \Sg_1$ and $(d)$ is proved.

The expression for the Hessian matrix of second derivatives in polar coordinates is
$$ D^2(\phi) = Q_\theta H(\phi) Q_\theta^T  
$$
where $Q_\theta$ is the rotation matrix by angle $\theta$, and
\begin{equation}
H(\phi) = \left(  \begin{array}{cc}
\frac{\partial^2 \phi}{\partial r^2}    & \frac{1}{r}\left(\frac{\partial^2 \phi}{\partial r \partial \theta}  
     - \frac{1}{r} \frac{\partial \phi}{\partial \theta} \right) \vspace{0,3cm}\\
   \frac{1}{r}\left(\frac{\partial^2 \phi}{\partial r \partial \theta}  
     - \frac{1}{r} \frac{\partial \phi}{\partial \theta} \right)  & 
   \frac{1}{r^2} \frac{\partial^2 \phi}{\partial \theta^2} + \frac{1}{r} \frac{\partial \phi}{\partial r} 
\end{array}
\right). \nonumber
\end{equation}

Now, we have that (recall that $x_0 = 0$)
$$   \phi(\theta, r ,t) = \frac{\bar{a}}{\rho(\theta)^2} r^2 - \beta t^2 +M.$$
One can notice that $\phi$ is well defined and smooth in $\Omo$, (which means $\{r \geq \E\} \setminus  \Gamma_1$). All the computations that follows are valid in this set. We already said above that $\rho$ is constant with respect to $r$ and only depends on $\theta$ such that $\frac{\partial \rho}{\partial r} = 0$. Hence, we have that 
\begin{equation}
\label{Hess}
H(\phi) = \frac{2 \bar{a}}{\rho^2}
\left(  \begin{array}{cc}
1   & -\frac{\rho_{\theta}}{\rho}  \\
-\frac{\rho_{\theta}}{\rho}   &  
\frac{1}{\rho^2}( 3\rho_{\theta}^2 - \rho\rho_{\theta \theta} + \rho^2)
\end{array}
\right),
\end{equation}
where we have denoted $\rho_{\theta} = \frac{\partial \rho}{\partial \theta}$ and so on.

We  will use the following well known facts (see for example \cite{Gray}) concerning curves in the plane:
\begin{lem} Let $\gamma$ be a $C^2$ curve in the plane. Then:
\begin{enumerate}
\item[(a)]  The curve $\gamma$ is strictly convex only at those points where their curvature is positive.
\item[(b)] If $\gamma$ is parameterized in polar coordinates by its angle, that is $$\gamma(\theta) = (r(\theta)\cos\theta, r(\theta)\sin\theta),$$ then the curvature of $\gamma$ is given by the formula
$$ \kappa_{\gamma}(\theta) = \frac{r^2 + 2r_{\theta}^2- r r_{\theta \theta}}{(r^2 + r_{\theta}^2)^{3/2}}.
$$
\end{enumerate}
\end{lem}
Since the polar parametrization of $\Gm_1$ is precisely the above function with $r(\theta) = \rho(\theta)$ and  $\Om_1$ is strictly convex, we obtain
 \begin{equation}
\label{curvat}
 \kappa_{_{\Gamma_1}} (\theta) = \frac{\rho^2 + 2\rho_{\theta}^2- \rho \rho_{\theta \theta}}{(\rho^2 + \rho_{\theta}^2)^{3/2}} > 0  \hspace{1cm} \forall \theta \in [0, 2 \pi[.
 \end{equation}
From (\ref{Hess}) and (\ref{curvat}) we have 
\begin{eqnarray*}
\Delta \phi & = & \textnormal{tr}(D^2(\phi)) \\
		& = & \textnormal{tr}(H(\phi)) \\
		& = & 2c \left( 1 + \frac{1}{\rho^2}( 3\rho_{\theta}^2 - \rho\rho_{\theta \theta} + \rho^2) 
		\right)\\
		& \geq & 2c \left( 1  + \frac{1}{\rho^2} (\rho^2 + \rho_{\theta}^2)^{3/2} \kappa_{_{\Gamma_1}} \right)~ > ~2c
\end{eqnarray*}
and $(e)$ is proved.

We also have
\begin{eqnarray*}
\det(H( \phi)) & = & \frac{2c}{\rho^2} \left(  3\rho_{\theta}^2 - \rho\rho_{\theta \theta} + \rho^2 - \rho_{\theta}^2 \right) \\
   & = & \frac{2c}{\rho^2} (\rho^2 + \rho_{\theta}^2)^{3/2}  \kappa_{_{\Gamma_1}}~ > ~ 0.
\end{eqnarray*}

By the Sylvester's Criterion we can see that $H(\phi)$ (thus $D^2(\phi)$) is positive definite. Indeed, the element $\frac{2a}{\rho^2}$ and the determinant of the matrix $H(\phi)$ are positive. Finally, since $\overline{\Om}$ is compact, this implies $(f)$ and the proof of Lemma~\ref{propphi} is complete.
\QED

We introduce now the last hypothesis we will need  in order to get the Carleman inequality:
\begin{equation} \label{H2} 
 \beta < \min\left\{ \frac{ \min\{a_1, a_2\} \delta_1}{2}, \frac{M}{T^2} \right\} 
\end{equation} 
 and we take $\gm \in (0,1)$  such that 
 \begin{equation} \label{H3} 
 \gm > \frac{2 \beta}{\beta + \frac{a_1a_2}{\textrm{diam}(\Omega)^2}}  
 \end{equation}
 \begin{equation} \label{H4} 
 \gm < \frac{ 2 \min\{a_1, a_2\}  \delta_1}{2 \beta + \max\{a_1, a_2\} \| \Delta \phi \|_{L^\infty(\Om_0)}^2}. 
  \end{equation}

\begin{rem}
\begin{enumerate}
\item We have to take $M$ large enough in order to $(\ref{H2})$ and the hypothesis of the inverse problem (see Theorem~\ref{teolp}) become compatible.
\item Taking $\beta$ small enough, $(\ref{H3})$ and $(\ref{H4})$ become compatible. But, as we will see in the next section, the smaller is $\beta$, the bigger need to be the inversion time for the inverse problem (see Theorem~\ref{teolp}).  
\item Actually, the optimal $\delta_1$ is the first eigenvalue of $D^2(\phi)$. Having an explicit expression for it could help to a better choice of $\beta$.
\item From the hypothesis of Theorem \ref{Carleman} we have $a_2 < a_1$ and the maximum and the minimum in $(\ref{H2})$ and $(\ref{H4})$ are known. 
\end{enumerate}
\end{rem}

\begin{rem} Pseudoconvexity.
It is easy to check that Proposition \ref{propphi} and assumption $(\ref{H2})$ imply that  $\phi_1$ and
and $\phi_2$ are pseudoconvex \cite{Albano-Tataru,Hor} with respect to $P$ in $\Om_1\setminus \Bo$ and $\Om_2$ respectively. Global Carleman estimates (without explicit dependence on the parameter $\lambda$) can be deduced in each subdomain $\Omega_1\setminus \Bo$ and $\Omega_2$, 
for solutions that vanish on the exterior boundary and the interface (see \cite{Tata}, \cite{Albano-Tataru}). Nevertheless, the traces of the transmission wave equation do not  vanish on the interface $\Gamma_1$ and we will need to use  carefully  hypothesis $a_1>a_2$ and the parameter $s$.
On the other hand, we have to use the parameter $\lambda$ in order to get rid of the lack of pseudoconvexity in $\Bo$.  
\end{rem}

\begin{rem} \label{Coefvar} Coefficient  $a = a(x)$ variable.
We can also prove a Carleman estimate in the more general case $a(x) = a_1 (x) \1_{\Om_1}  + a_2(x)\1_{\Om_2}$ with $a_j \in C^1(\overline \Om_j)$, $j=1,2$, if each $a_j$ is constant in the interface $\Gamma_1$, and under some hypothesis on $\nabla a$ (similar to those of \cite{ImYam1}). More precisely, if we check the pseudoconvexity condition in this case, we will have that $\phi$ is pseudoconvex with respect to the operator $\partial_{tt} - div(a\nabla u)$ (in each domain $Q_j$) if there exists $\theta \in (0,1)$ such that:
\begin{eqnarray}
\frac{|\nabla \phi \cdot \nabla a_j|}{2a_j}  <  \delta_1(1- \theta),\ j=1,2 \label{cond-a}\\
\beta  < \frac{\sqrt a_j}{\sqrt a_j + T |\nabla a_j|} \left( \frac{a_j \delta_1 \theta}{2} \right),\ j=1,2.
\label{cond-beta}
\end{eqnarray}
It is easy to check that these hypothesis are compatible with the assumptions of Theorem~\ref{ip}.
Indeed, for $T$ sufficiently large, there exists $\beta$ satisfying both $\beta>a_1/T^2$ and (\ref{cond-beta}).  Nevertheless,  in order to construct the weight function as we have done above (and deal with the traces of the solutions on the interface $\Gamma_1$), it is crucial for each function $a_j$ to 
be constant on the interface.
\end{rem}

\subsection{Listing all the terms}\label{list}
In this part of the work we develop the $L^2$-product of $P_1(w)$ and $P_2(w)$. We will do formal computations, by writing generically $\phi$ for the weight function and  $Q$ for the domain with boundary $\Sigma$. 

As presented in Section 1, we have, for $\lambda > 0$, $s > 0$,
$$\varphi  =  e^{\lambda \phi} ~,~~
w =  e^{s \varphi} u$$
$$P(w)  =   e^{s \varphi} L( e^{- s \varphi} w)
          =  P_1(w) + P_2(w) + R(w) $$
where
\begin{eqnarray*}
P_1(w) & = & w_{tt} - a\Delta w + s^2\lambda^2\varphi^2E(\phi)w, \\
P_2(w) & = & (\gm-1)s \lambda \varphi L(\phi) w - s \lambda^2 \varphi E(\phi) w, 
- 2s \lambda \varphi (\phi_t w_t - a \nabla \phi \cdot \nabla w) .\\
R(w)  & =  & - \gm s \lambda \varphi L(\phi)w 
\end{eqnarray*}

We set  $\left<P_1(w),P_2(w)\right>_{L^2} = \sum\limits_{i,j=1}^3 I_{i,j}$, where  $I_{i,j}$ is the integral of the product of the $i$th-term in $P_1(w)	$ and the $j$th-term in $P_2(w)$. Therefore,
\begin{eqnarray*}
         I_{1,1}   & = &  -s \la (\gm -1) \iint_Q  \varphi L(\phi) |w_t|^2  + \frac{s\la^2(\gm-1)}{2} \iint_Q  |w|^2 \varphi (\phi_{tt} + \la |\phi_t|^2)L(\phi), \\
\smallskip
I_{1,2} & = & s\la^2 \iint_Q  |w_t|^2 \varphi E(\phi) - s \la^2 \iint_Q  |w|^2 \varphi |\phi_{tt}|^2 
              - \frac{5s\la^3}{2} \iint_Q  |w|^2 \varphi \phi_{tt} |\phi_t|^2 \\
    	&   & + \frac{s \la^3}{2} \iint_Q  |w|^2a \varphi \phi_{tt}|\nabla \phi|^2  -  \frac{s\la^4}{2}\iint_Q  |w|^2 \varphi |\phi_t|^2 E(\phi), 
\end{eqnarray*}
\begin{eqnarray*}
I_{1,3} & = & s\la \iint_Q  |w_t|^2\varphi (\phi_{tt} + \la |\phi_t|^2)  - 2 s\la^2 \iint_Q  w_t \phi_t \varphi a \nabla w \cdot \nabla \phi \\
	&   & + s\la \iint_Q  |w_t|^2 \varphi a ( \Delta \phi + \la |\nabla \phi|^2)  -  s \la \iint_{\Sigma} |w_t|^2a \varphi \nabla \phi \cdot \nu , \\
\smallskip
I_{2,1} & = & - s \la (\gm -1) \iint_{\Sigma} \varphi L(\phi) w a \dn{w}  + s \la (\gm -1) \iint_Q a |\nabla w|^2 \varphi L(\phi) \\
	&   & - s\la^2 \frac{\gm -1}{2}\iint_Q |w|^2 \varphi a L(\phi)(\la |\nabla \phi|^2 + \Delta \phi) 
              - s\la \frac {\gm -1}{2} \iint_Q |w|^2 a \varphi \Delta(L(\phi))  \\
	&   & - s\la^2 (\gm -1) \iint_Q |w|^2 a \varphi (\nabla \phi \cdot \nabla L(\phi) )
 		+ s\la \frac {\gm -1}{2} \iint_{\Sigma} |w|^2a \varphi \nabla L(\phi) \cdot \nu \\
	&  &	+ s\la^2 \frac{\gm -1}{2} \iint_{\Sigma}|w|^2a \varphi L(\phi) \dn{\phi},
\end{eqnarray*}
\begin{eqnarray*}
I_{2,2} & = & s\la^2 \iint_{\Sigma} a \varphi E(\phi) w \nabla w \cdot \nu  
		- \frac{s\la^2}{2} \iint_{\Sigma} a |w|^2 ( \la \varphi E(\phi) \nabla\phi + \varphi \nabla E(\phi))\cdot \nu \\
	&   &	+  \frac{s\la^3}{2} \iint_Q |w|^2 a \varphi E(\phi)(\Delta \phi + \la|\nabla \phi|^2) 
		- 2s \la^3 \iint_Q |w|^2 a^2 \varphi D^2(\phi)(\nabla \phi, \nabla \phi)\\
	&   &	+ \frac{s\la^2}{2} \iint_Q |w|^2 a \varphi \Delta( E(\phi))  -  s\la^2 \iint_Q |\nabla w|^2 a \varphi E(\phi) ,
\end{eqnarray*}
 \begin{eqnarray*}
I_{2,3} & = & s\la \iint_Q |\nabla w|^2 a \varphi L(\phi) + s\la^2\iint_Q |\nabla w|^2 a \varphi E(\phi) 
		+ 2s\la^2 \iint_Q a^2 \varphi |\nabla \phi \cdot \nabla w|^2 \\
	&   &  - 2s\la^2 \iint_Q a \varphi \phi_t w_t \nabla w\cdot\nabla\phi + 2s\la \iint_Q a^2 \varphi D^2(\phi)(\nabla w, \nabla w)\\
	&   &  + s\la \iint_{\Sigma}|\nabla w|^2a^2 \varphi \nabla \phi \cdot \nu
		+ 2s\la \iint_{\Sigma} a \varphi (\phi_t w_t - a \nabla \phi \cdot \nabla w ) \dn{w}, \\
\smallskip
I_{3,1} & = &  s^3 \la^3 (\gm -1) \iint_Q |w|^2 \varphi^3 L(\phi) E(\phi), \\
\smallskip
I_{3,2} & = & - s^3 \la^4  \iint_Q |w|^2 \varphi^3 E(\phi)^2, \\
\smallskip
I_{3,3} & = &  s^3 \la^3  \iint_Q |w|^2 \varphi^3 E(\phi)L(\phi) + 2 s^3 \la^3  \iint_Q |w|^2 \varphi^3(|\phi_t|^2\phi_{tt} +a^2D^2(\phi)(\nabla \phi, \nabla \phi) )  \\
	&   & + 3 s^3 \la^4  \iint_Q |w|^2 \varphi^3 E(\phi)^2 +
   s^3\la^3 \iint_{\Sigma} a |w|^2 \varphi^3E(\phi) \dn{\phi}.
\end{eqnarray*}
Gathering all these terms,we get
\begin{eqnarray*}
\left<P_1(w), P_2(w)\right>_{L^2(Q)} & = & 2s\la \iint_Q |w_t|^2 \varphi \phi_{tt} 
			     -\gm s\la \iint_Q |w_t|^2  \varphi L(\phi) \\
			    &    & +2s \la^2 \iint_Q \varphi \left( |w_t|^2 |\phi_t|^2  - 2 w_t \phi_t a \nabla w \cdot \nabla \phi + a^2 |\nabla \phi \cdot \nabla w|^2 \right) \\
			    &    & +2s \la \iint_Q a^2 \varphi D^2(\phi)(\nabla w, \nabla w) \\
			    &    &  + \gm s \la \iint_Q a |\nabla w|^2 \varphi L(\phi) 
			       + 2s^3\la^4   \iint_Q |w|^2 \varphi^3 E(\phi)^2 \\
			    &    &  + 2s^3\la^3  \iint_Q |w|^2 \varphi^3	 (|\phi_t|^2\phi_{tt} + a^2 D^2(\phi)(\nabla \phi, \nabla \phi) )\\		
 			    &    &  + \gm s^3\la^3  \iint_Q |w|^2 \varphi^3L(\phi)E(\phi)\\
			    &	   & + X + J,
\end{eqnarray*}
where $J$ is the sum of all the boundary terms:
\begin{eqnarray}
J & = & s \la \inS \left( a^2 \varphi |\nabla w |^2 \dn{\phi}  
			    - 2 a^2 \varphi (\nabla \phi \cdot \nabla w) \dn{w} \right)\nonumber \\
  &   & + s \la \frac{\gm-1}{2} \inS |w|^2\varphi a \nabla L(\phi)\cdot \nu \nonumber  \\
 &   & + s \la^2 \frac{\gm-1}{2} \inS |w|^2\varphi a  L(\phi)  \dn{\phi} \nonumber \\
  &    & - s \la (\gm-1) \inS w a \dn{w} \varphi  L(\phi) \nonumber 
  + s \la^2 \inS a \varphi E(\phi) w \dn{w}  \label{interface} \\
  &  &  - s \la^3 \frac{1}{2} \inS |w|^2  \varphi E(\phi) a  \dn{\phi} \nonumber \\
 & &   - s \la^2 \frac{1}{2} \inS |w|^2 a \varphi \nabla E(\phi)\cdot \nu \nonumber \\
  &  &  + 2 s \la \inS a \varphi \phi_t w_t \dn{w} \nonumber 
   -  s \la\inS |w_t|^2\varphi a \dn{\phi} \nonumber \\
  &  & + s^3 \la^3 \inS |w|^2 \varphi^3 E(\phi) a \dn{\phi} \nonumber 
  \end{eqnarray}
and $X$  is a the sum of the remaining terms, in such a way that:
$$
|X| \leq Cs\la^3 \iint_Q \varphi^3|w|^2
$$

In the sequel, we denote by $A_j$, $j = 1,..., 8$ the first eight integrals we have 
listed in the product of $P_1(w)$ by $P_2(w)$. Thus, we have
\begin{equation}\label{prodL2}
 \left<P_1(w), P_2(w) \right>_{L^2(Q)} = \sum\limits_{j=1}^8 A_j  + X + J.
 \end{equation}

 \subsection{Proof of Theorem \ref{Carleman}} \label{proofCarlem}
We take $\beta$, $\gamma$ and $M$ satisfying the hypothesis $(\ref{H1}), (\ref{H2}), (\ref{H3})$ and $(\ref{H4})$, and $\phi$ as the corresponding weight function. We assume throughout  all this part of the work the hypothesis of 
Proposition \ref{propphi} (especially that $\Om_1$ is strictly convex with $C^3$ boundary).

Recall the notation $Q_0 = (\Om_1 \cup \Om_2) \times (-T,T)$ and $Q_{x_0} = (\Om_0\setminus \Bo) \times (-T,T)$. We apply the above computations to $w = e^{s \varphi}u$ with $u \in X$ in each one of the open sets $Q_1$ and $Q_2$. Adding the terms that result in both cases, we have (recall that $u(\pm T) = u_t(\pm T) = 0$ for all $u \in X$ ):
 \begin{equation}\label{prodL2Q}
 \left<P_1(w), P_2(w) \right>_{L^2(Q_0)} = \sum\limits_{j=1}^8 A_{j,Q_0}(w)  + X_{Q_0}(w) + J_{\Sigma_1}(w) + J_{\Sigma}(w),
 \end{equation}
where we have written $A_{j,Q_0}$ instead of the integral $A_j$ given in subsection \ref{list} taken in the set $Q_0$, etcetera.
  
The proof of Theorem \ref{Carleman} is based on the next three facts:
 \begin{itemize}
 \item The sum of the $A_j$-integrals in $Q_{x_0}$ can be minored.
 \item The sum of terms in the interface given by $J_{\Sigma_1}$ is nonnegative, and:
 \item We can introduce a second weight function centered at a point different to $x_0$ in order to deal with the integrals in $\Bo$. 
 \end{itemize}
 
 The key points in each step of the proof are based on the properties of $\phi$ listed in Proposition \ref{propphi}.
 \subsubsection{The interior}
 \begin{propo} \label{propinterior}
 There exist $\delta > 0$, $C>0$ and $\la_0 >0$ such that:
 \begin{equation*}
 \sum_{j=1}^8 A_{j, Q_0}(w)  \geq  \delta \| w \|_{\Om_0, \varphi}  - C \|w\|_{\Bo, \varphi}
  \end{equation*}
for all $\la \geq \la_0$, for all $u \in X$.
 \end{propo}
 
 \dem 
 We arrange the terms into four groups: 
 \begin{enumerate}
 \item $ A_{1,Q_0} + A_{2,Q_0}  = s \la \ds\iint_{Q_0} |w_t|^2 \varphi (-\gm L(\phi) - 4 \beta)$.
For all $(x,t)$ in $Q_{x_0}$ we have: 
 \begin{eqnarray*}
 -\gm L(\phi) - 4 \beta & = & \gm (2 \beta +a \Delta \phi) - 4 \beta\\
 				  & \geq & \gm (2 \beta + 2a c(x)) - 4 \beta  \hspace{0.9cm} \mbox{by Proposition \ref{propphi} }\\
				   & \geq & \gm (2 \beta + 2a \frac{\bar{a}}{\textrm{diam}(\Omega)^2}) - 4 \beta  \hspace{0.7cm} \mbox{  by definition of $c(x)$ }\\
				  & =  & \delta > 0 \hspace{3.15cm}  \mbox{   by $(\ref{H3})$.}
 \end{eqnarray*}
 Therefore:
 $$
 A_{1,Q_0} + A_{2,Q_0}  \geq \delta_1  s \la \iint_{\Omo} |w_t|^2 \varphi   
-  C  s \la \iint_{\Bo} |w_t|^2 \varphi  $$
 
 \item $A_{3,Q_0}
  = 2s\la^2 \ds\iint_{Q_0} \varphi (\phi_t w_t - a \nabla \phi \cdot \nabla w)^2 \geq 0 $
 
 \item $A_{4,Q_{x_0}} + A_{5,Q_{x_0}} = s \la \ds\iint_{Q_{x_0}} \varphi \left(a^2 D^2(\phi)(\nabla w, \nabla w) + \gm a L(\phi)|\nabla w|^2 \right)$.
Then, using Proposition \ref{propphi} and  $(\ref{H4})$, we obtain
\begin{eqnarray*}
  A_{4,Q_{x_0}} + A_{5,Q_{x_0}}   &  \geq & s \la \iint_{Q_{x_0}} \varphi \left(2a^2 \delta_1|\nabla w|^2  + \gm a L(\phi)|\nabla w|^2 \right)  \\ 
   &  \geq  &  s \la \iint_{Q_{x_0}} \varphi a\left( 2a \delta_1 - \gm(2 \beta + a \Delta \phi) \right)  |\nabla w|^2 \\
   &  \geq  &  s \la \iint_{Q_{x_0}} \varphi a \left( 2a \delta_1 - \gm(2 \beta + a \| \Delta \phi\|_{L^{\infty}}) \right)  |\nabla w|^2 \\
   &  \geq  &  s \la \delta_2 \iint_{Q_{x_0}} \varphi   |\nabla w|^2
    \end{eqnarray*}
Therefore
 $$
 A_{4,Q_{0}} + A_{5,Q_{0}}  \geq  \delta_2  s \la \iint_{\Omo} \varphi   |\nabla w|^2 
 - C   s \la \iint_{\Bo} \varphi   |\nabla w|^2   
 $$
 
 \item $ \sum\limits_{j=6}^8 A_{j,Q_0}    =  s^3 \la^3 \ds\iint_{Q_0} |w|^2 \varphi^3 F_{\la}(\phi) ~$ where
 \begin{eqnarray*}	
F_{\la} (\phi)& = & 2 \la E(\phi)^2 +2|\phi_t|^2\phi_{tt} + 2a^2 D^2(\phi)(\nabla \phi, \nabla \phi) + \gm L(\phi)E(\phi)\\
 &= & 2 \la E(\phi)^2 + \gm L(\phi)E(\phi) - 16 \beta^3 t^2 + 2a^2 D^2(\phi)(\nabla \phi, \nabla \phi) \\
           & = &2 \la E(\phi)^2 + \underbrace{(\gm L(\phi)-4\beta )}_{b(x) < 0} E(\phi) - 4 \beta a|\nabla \phi|^2 +  2a^2 D^2(\phi)(\nabla \phi, \nabla \phi)
            \end{eqnarray*}
From Proposition \ref{propphi} and  $(\ref{H2})$, there exists $d_0 > 0$ such that for all
$(x,t) \in Q_{x_0}$ we have
 \begin{eqnarray*}	           
 	F_{\la} (\phi)& \geq & 2 \la E(\phi)^2 + b(x)E(\phi)  
			- 4 \beta a|\nabla \phi|^2 + 2a^2 \delta_1|\nabla \phi|^2 \\
 	& \geq & 2 \la E(\phi)^2 + b(x)E(\phi)  + a(2a \delta_1 - 4 \beta ) |\nabla \phi|^2  \\
 	& \geq & 2 \la E(\phi)^2 - \|b\|_{\infty} |E(\phi)|  +  a(2a \delta_1 - 4 \beta ) |\nabla \phi|^2 \\
 	& \geq & 2 \la E(\phi)^2 - \|b\|_{\infty} |E(\phi)|  + d_0 \\
	& = & \left\{ \begin{array}{rr} f_{1,\la}(E(\phi)) & \mbox{  if } E(\phi) > 0 \\
						f_{2,\la}(E(\phi)) & \mbox{ if } E(\phi) < 0
			  \end{array}  \right.
 \end{eqnarray*}
where  \begin{eqnarray*}
f_{j,\la} : \RR &\longrightarrow& \RR \\
x & \longmapsto & 2\la x^2  + (-1)^j \|b\|_{\infty} x +d_0.
\end{eqnarray*}
As $d_0 > 0$, there exists $\la_0>0$ such that for all $\la \geq \la_0$ 
 $$
  \min_{\RR} (f_{j,\la})  \geq  \frac{d_0}{2}  > 0  \hspace{0.6cm} j =1,2.
 $$
Thus, for each $\la \geq \la_0$ we have
 $$
  \sum\limits_{j=6}^8 A_{j,Q_0} 
      \geq \delta_3 s^3 \la^3 \iint_{Q_{x_0}} |w|^2 \varphi^3  - C s^3 \la^3  \iint_{\Bo} |w|^2 \varphi^3.
 $$
 
 \end{enumerate}
By collecting all the terms $A_{j,Q_0}$ together, we conclude the proof of Proposition \ref{propinterior}.
\QED

 \subsubsection{The interface} 
 Since the interface $\Gamma_1$ is a common boundary of $\Om_1$ and $\Om_2$, the term 
 $J_{\Sigma_1}$ in $(\ref{prodL2Q})$ is the sum of the integrals comming from each domain:
 $J_{\Sigma_1} = J_{\Sigma_1}(w_1) + J_{\Sigma_1}(w_2)$. We have the following result:
 \begin{propo} \label{propinterface}
 Suppose $0 < a_2 < a_1$. Then there exists $s_0 >0$ such that 
 $$
J_{\Sigma_1} = J_{\Sigma_1}(w_1) + J_{\Sigma_1}(w_2)    \geq 0  \hspace{1cm} \forall s \geq s_0
$$
for all $u \in X$.
\end{propo}
 
 \dem
We enumerate the ten integrals arising in (\ref{interface}) associated with the common  boundary $\Sigma_1$, and we denote by $J_i$ the sum of the $i$-th integral in (\ref{interface}) which comes from $\Om_1$ with the respective one of $\Om_2$: $J_i = J_i(w_1) + J_i(w_2)$.
 
 In order to prove the inequality, we arrange the terms into three groups. In each case, we use Proposition \ref{propphi} and the fact that $w$ satisfies the transmission conditions.
 \begin{enumerate}
 
 \item  Is not difficult to see that $  J_k = 0 $ for each $k \in \{3, 4, 7, 8, 9\}$. Indeed, from (d) of Proposition \ref{propphi} we get $L_1(\phi_1) =  L_2(\phi_2)$ and $a_1\nabla E_1(\phi_1) = a_2\nabla E_2(\phi_2)$ on $\Sigma_1$, and the desired result follows.
 
Now, let us denote by $g$ the real function defined in $\Sigma_1$ by
 $$
 g(x,t) := E_1(\phi_1) - E_2(\phi_2) = \left(\frac{1}{a_2} - \frac{1}{a_1}\right)\left|a_1\dns{\phi}{1}\right|^2
 $$
 
 Since $a_2 < a_1$, we have $g > 0$ in $\Sigma_1$.
 
 Thus we can prove:
 
 \item $J_2 + J_6 + \frac{1}{2}J_{10} \geq 0 \hspace{1.5cm}  \forall s \geq s_0 $ \\
Indeed:
\begin{eqnarray*} 
 -J_2 - J_6 & =   &  s\lambda \frac{1-\gm}{2} \iint_{\Sg_1}|w|^2 \varphi a \nabla L(\phi)\cdot \nu \\
  		&    & + s \la^3 \frac{1}{2} \iint_{\Sigma_1} |w|^2  \varphi \left( a_1  \dns{\phi_1}{1}\right) g(x,t)  \\
 	& \leq & \frac{1}{2}J_{10}
  \end{eqnarray*}
   for all $s \geq s_0$, since $\varphi \geq 1$.
 
 \item $ J_1 + J_5 + \frac{1}{2}J_{10} \geq 0 \hspace{1.5cm}  \forall s \geq s_0 $
 
 By construction, $\phi$ is constant on each 
 level $\Gm_1 \times \{t\}$ of the interface (Proposition \ref{propphi}). Thus
 \begin{equation}\label{grad}
 \nabla \phi_j \cdot \nabla w_j = \dns{\phi_j}{j} \dns{w_j}{j} \quad \mbox{  in  } \Sg_1 
 \, \, \mbox{ for } j = 1,2.
 \end{equation}
 Moreover, since $w$ satisfies (\ref{Tr}) we have
 $$
 \left|\frac{\partial w_1}{\partial \tau_1} \right| = \left|\frac{\partial w_2}{\partial \tau_2} \right | \quad \mbox{ in } \Sg_1 \, \, \forall  u \in X. 
 $$
 
Hence:
\begin{equation}\label{tan}
\sum\limits_{j=1}^2
\left |\frac{\partial w_j}{\partial \tau_j} \right |^2 a_j^2\varphi_j \dns{\phi_j}{j}
= \left |\frac{\partial w_1}{\partial \tau_1} \right |^2 \varphi_1\left(a_1 \dns{\phi_1}{1} \right)( a_1 - a_2)  > 0
\end{equation}

From (\ref{grad}) and (\ref{tan}) we get:
\begin{eqnarray*}
J_1 & \geq & - s \la \iint_{\Sg_1}  \left( \left|a_1\dns{w_1}{1}\right|^2 \varphi_1\dns{\phi_1}{1} + \left|a_2\dns{w_2}{2}\right|^2 \varphi_2\dns{\phi_2}{2} \right) \\
    & = &  - s \la \iint_{\Sg_1} \left|a_1\dns{w_1}{1}\right|^2 \varphi_1 \left( \dns{\phi_1}{1} + \dns{\phi_2}{2} \right) \\
    & = &  s \la \iint_{\Sg_1} \left|a_1\dns{w_1}{1}\right|^2 \varphi_1 \left( \frac{1}{a_2} - \frac{1}{a_1} \right) \left( a_1 \dns{\phi_1}{1} \right).   
\end{eqnarray*}

On the other hand,
 \begin{eqnarray*}
 -J_5 & = & - s \la^2 \iint_{\Sg_1}  w_1 \varphi_1 \left(a_1 \dns{w_1}{1} \right)g(x,t) \\
 	& \leq & \frac{1}{2}s^2 \la^3 \iint_{\Sg_1}  |w_1|^2 \varphi_1 g + 
	\frac{1}{2} \la  \iint_{\Sg_1}  \varphi_1 g \left|a_1 \dns{w_1}{1} \right|^2 \\
	& = & \frac{1}{2}s^2 \la^3 \iint_{\Sg_1}  |w_1|^2 \varphi_1 g \\
	&    & + \frac{1}{2} \la  \iint_{\Sg_1} \left|a_1 \dns{w_1}{1} \right|^2 \varphi_1 \left(\frac{1}{a_2} - \frac{1}{a_1}\right)\left|a_1\dns{\phi}{1}\right|^2\\
&  \leq & \frac{1}{2 \delta} s^2\la^3 \iint_{\Sg_1} |w|^2  \varphi_1 g \left(a_1 \dns{\phi_1}{1}\right)  \\
	&  & + \frac{C}{2} \la \iint_{\Sg_1} \left|a_1 \dns{w_1}{1} \right|^2 \varphi_1 \left(\frac{1}{a_2} - \frac{1}{a_1}\right) \left(a_1\dns{\phi}{1}\right)\\
	& \leq  &  \frac{1}{2}J_{10}  +  J_1  \quad \forall s \geq s_1.
\end{eqnarray*}
	
 \end{enumerate}
 Proposition (\ref{propinterface}) is proved.
 \QED
 
 \subsubsection{The boundary $\Sigma$.} 
  Since we deal with functions $w$  such that $w_2 = 0$ in $\Sigma$, we have  
\begin{eqnarray} 
 J_{\Sigma}   = J_1(w_2) & = &-s \la \inS \varphi \left|a_2\dn{w} \right|^2 \left(\dn{\phi}\right) \nonumber \\
	& \geq & -s \la \iint_{\Sigma_+}\varphi \left|a_2\dn{w} \right|^2 \left(\dn{\phi}\right) \nonumber \\
	& \geq & -s \la \iint_{\Sigma_+} \varphi \left|a_2\dn{w} \right|^2 \left\|\dn{\phi}\right\|_{L^{\infty}(\Sg)} 
			\label{bd}\\
 	&  =   & -s \la C\iint_{\Sigma_+} \varphi \left|a_2\dn{w} \right|^2, \nonumber
 \end{eqnarray}					   
where we have defined $\Sigma_+ = \{ (x,t) \in \Gamma \, : \nabla \phi(x,t) \cdot \nu(x) > 0 \}$.

 \subsubsection{Carrying all together.} 
From (\ref{prodL2Q}), (\ref{bd}) and Propositions \ref{propinterior} and 
\ref{propinterface}, there exist $s_0$, $\la_0$, $C \in \RR$ 
such that for each $s \geq s_0$ and $\la \geq \la_0$ we have
\begin{eqnarray}
 \left\|w\right\|_{_{\Om_0,\varphi}}^2   - C \left\|w\right\|_{_{\Bo,\varphi}}^2 
& + X_{Q_0}  \nonumber   \\
  - s \la C \iint_{\Sigma_+} \varphi \left |a_2\dn{w} \right |^2   & \leq \quad  C \left< P_1(w), P_2(w) \right>_{L^2(Q_0)}  \label{ineg1}
\end{eqnarray} 

Adding $ \frac{C}{2} \left( |P_1(w)|^2_{L^2(Q_0)} + |P_2(w)|^2_{L^2(Q_0)}\right) $ at 
both sides of $(\ref{ineg1})$ we obtain
\begin{eqnarray}
  \left|P_1(w) \right |^2_{L^2(Q_0)} + \left|P_2(w) \right|^2_{L^2(Q_0)}
 +   \left\|w\right\|_{_{\Om_0,\varphi}}^2   - C & \left\|w \right\|_{_{\Bo,\varphi}}^2 \nonumber \\
    + X_{Q_0}  - s \la C \iint_{\Sigma_+} \varphi \left |a_2\dn{w} \right |^2   
       &  \leq  C \left |P_1(w) + P_2(w) \right |^2_{L^2(Q_0)}  \nonumber \\
       &  =  C \left |P(w) - R(w) \right |^2_{L^2(Q_0)}  \label{ineg2} \\ 
       &   \leq  C \left(\left |P(w) \right |^2_{L^2(Q_0)}    
      			+ \left |R(w) \right|^2_{L^2(Q_0)}\right) \nonumber
 \end{eqnarray} 
Thanks to $(\ref{H2})$ we have $\lambda \leq C \varphi$ for $\lambda$ large enough. 
Therefore
\begin{eqnarray} 
|X_{Q_0}| + C|R(w) |^2_{L^2(Q)} & \leq  C s \la^3 \iint_{Q} \varphi^2 |w|^2 
					\quad & \forall \la \geq \la_1 \nonumber \\ 
	&	 \leq  \frac{1}{2}\left \|w\right\|_{_{\Om_0, \varphi}}^2	\quad & \forall s \geq s_1.
		 \label{ineg3} 
\end{eqnarray}

From (\ref{ineg2}) and (\ref{ineg3}) we get, for all $s \geq \max \{s_0, s_1\}$, 
$\la \geq \max \{\la_0, \la_1\}$:
\begin{eqnarray}
  |P_1(w)|^2_{L^2(Q_0)} + |P_2(w)|^2_{L^2(Q_0)}
+ \left \|w\right\|_{_{\Om_0,\varphi}}^2    & \leq C |P(w)|^2_{L^2(Q_0)}  \nonumber \\
& + C \left\|w\right\|^2_{_{\Bo, \varphi}} 
 +   s \la C  \iint_{\Sigma_+} \varphi \left|a_2\dn{w}\right |^2.  \label{Carl0}
\end{eqnarray}

\subsubsection{Eliminating the term in $B_{\E}(x_0)$.}
In the last step we will remove the integral in $B_{\E}(x_0)$ from the right hand side of (\ref{Carl0}). In order 
to do that, first remark that $x_0$ can be arbitrarily chosen in $\Om_1$ since $\Om_1$ is strictly convex.

Thus, we can take two different points in $\Om_1$ and we have the two respective inequalities given by (\ref{Carl0}). Now, we will show that the left hand side of each inequality can absorb the
term $\| \cdot \|_{\Bo}$ from the other inequality provided that $\E$ is small and $\la$ is large enough:

Denote by $x_1$, $x_2$ two points in $\Om_1$, and $\phi^1$, $\phi^2$ their respective weight functions. 
In order to have $\|\cdot \|_{\Bi,\varphi^1}$ absorbed by the term $ \|w^2 \|_{_{\Om_0,\varphi^2}}$ it 
suffices that
$$
C \varphi^1 < \frac{1}{2} \varphi^2  \hspace{1cm}  \mbox{in} \quad \Bi
$$
i.e.
$$
e^{\la(\phi^2 - \phi^1)} > 2C \hspace{1cm}  \mbox{in}  \quad \Bi
$$

Thus, if we show that it is possible to have $\phi^2 - \phi^1 > \delta > 0$ in $\Bi$ by 
taking $\la$ large enough we are done. 

In fact, let be $d = \frac{1}{2}|x_1 - x_2|$ and assume that $\E < d$. Then, for all $x \in \Bi$ we have:
\begin{eqnarray}
\phi^1(x,t) & \leq & \frac{a}{\rho_1^2} \E^2 - \beta t^2 + M \nonumber\\
		& \leq & \frac{a}{\alpha_1^2} \E^2 - \beta t^2 + M, \label{inegalpha}
\end{eqnarray}
where $\alpha_1 = d(x_1,\Gm_1) > 0$.

In the same way, if we denote $D_2 = \max\limits_{y \in \Gm_1}d(y,x_2)$, we get

\begin{eqnarray*}
\phi^2(x,t) & \geq &\frac{a}{\rho_2^2} d^2 - \beta t^2 + M \\
		&   \geq & \frac{a}{D_2^2}d^2 - \beta t^2 + M 
\hspace{1cm} \forall x \in \Bi.
\end{eqnarray*}
Consequently, we have
 \begin{equation}
\phi^2 - \phi^1  \geq  a \left(  \frac{d^2}{D_2^2} -  \frac{\E^2 }{\alpha_1^2} \right) 
\hspace{0.5cm} \forall x \in B_{\E}(x_1).
\end{equation}
It is clear that an analogous result is true by interchanging $x_1$ and $x_2$ (now with $\alpha_2$ and $D_1$). 
Thus, taking $\E < \min\left( \frac{d \alpha_1}{D_2}, \frac{d \alpha_2}{D_1} \right) $ we can absorb the desired  terms in the inequality and Theorem~\ref{Carleman} is proved.
\QED

\section{Proof of the stability of the inverse problem}\label{PbI}
In this section we apply the Carleman inequality of Theorem \ref{ip} to the inverse problem presented 
in Section \ref{Intro}. For a principal coefficient $a$ piecewise constant 
and $p\in L^\infty(\Omega)$, we consider the wave equation
\begin{equation}\label{W}
 \left\{ \begin{array}{rclc}
         u_{tt}  - \textnormal{div}( a(x) \nabla u) + p(x)u & = & g(x,t) & \Omega \times (0,T)\\
                              u & = & h & \Gamma \times (0,T)\\
                           u(0) & = & u_0  & \Omega\\
                           u_t(0) & = & u_1 & \Om.
\end{array}
\right. 
\end{equation}
If $g\in L^{1}(0,T;L^{2}(\Omega))$, 
$h\in L^{2}(0,T;L^{2}(\Gamma))$ and  $u_0\in H^1_0(\Omega)$, $u_1\in L^2(\Omega)$, then 
\cite{L-M,Da-Li} equation (\ref{W}) has a unique weak solution $u \in C([0,T];H^1_0(\Omega))\cap C^1([0,T];L^2(\Omega)) $ with continuous dependence in initial conditions and such that $\frac{\partial u}{\partial\nu}\in L^{2}(0,T;L^{2}(\Gamma))$.

In order to prove the local stability of the nonlinear application (\ref{nonlinear}), that is, the
problem of determining the potential $p$ in $\Omega$ by a single measurement of the flux $a_2 \dn{u}$ on $\Gamma$ between $t=0$ and $t=T$, we follow the ideas of \cite{Bu-Kli} and \cite{Kli-Mal}
Thus, we will first consider a linearized version of this problem, what means working on the wave equation
\begin{equation}
\label{lp}
 \left\{ \begin{array}{rclc}
         y_{tt}  - \textnormal{div}( a(x) \nabla y) + p(x)y  & = & f(x)R(x,t) & \Omega \times (0,T)\\
                              y & = & 0 & \Gamma \times (0,T)\\
                           y(0) & = & 0  & \Omega\\
                           y_t(0) & = & 0 & \Om
\end{array}
\right. 
\end{equation}
given $p$ and $R$, and proving the stability of the application $f|_\Omega\longmapsto \left.\frac{\partial y}{\partial \nu}\right|_\Sigma$. 

We will indeed prove the following result:
\begin{teo} \label{teolp}
For $x_1, x_2 \in \Om_1$ let  $R_j = \sup\{ |x-y_j(x)| \, : \, x \in \Om_2 \}$, $j=1,2$, where $y_j$ is 
defined in (\ref{defy}) with $x_0 = x_j$. Set $\alpha_j = d(x_j, \Gamma_1)$
and $\ds{ D_0 = \max \left\{ \frac{R_1+\alpha_1}{\alpha_1} , \frac{R_2 + \alpha_2}{\alpha_2} \right\} }$.  With the hypothesis of Theorem \ref{ip},  and $T$, $\beta$ satisfying $(\ref{H1}), (\ref{H2}) , (\ref{H3}) $ and $(\ref{H4})$,  suppose that
\begin{itemize}
\item $\|p\|_{L^\infty(\Om)} \leq m$
\item $T > D_0 \sqrt \frac{a_1}{\beta}$
\item $R \in H^1(0,T;L^{\infty}(\Om))$
\item $0 < r < |R(x,0)|$ almost everywhere in $ \Om$.
\end{itemize}
Then there exists $C>0$ such that for all $~ f \in L^2(\Om)$, the solution $y$ of 
$(\ref{lp})$ satisfies
$$  \|f\|^2_{L^2(\Om)} \leq C ~\left\| a_2 \dn{y} \right\|^2_{H^1(0,T;L^2(\Gm))}.  $$
\end{teo}

\dem
For each $f \in L^2(\Om)$ and  $R \in H^1(0,T;L^{\infty}(\Om))$, let $y$ be the solution of (\ref{lp}). We
take the even extension of $R$ and $y$ 
 to  the interval $(-T,T)$. We call this functions in the same way, and in this proof we denote $Q = \Om \times (-T,T)$ and $\Sg = \Gm \times (-T,T)$ the extended domains. Therefore, $z = y_t$ satisfies the following equation:
\begin{equation}
 \left\{ \begin{array}{rclc}
        z_{tt}  - \textnormal{div}( a \nabla z) + pz  & = & f(x)R_t(x,t) & Q\\
                              z & = & 0 & \Sg\\
            z(0) &=&0& \Omega\\
             z_t(0)  &=&  f(x)R(x,0)  & \Omega
\end{array}
\right. 
\end{equation}
and we have the usual energy estimate
$$
\|z\|_{H^1(-T,T;L^2(\Om))} \leq C \|fR_t\|_{L^2(-T,T;L^2(\Om))} + C \|fR(0)\|_{L^2(\Om)}
$$
that gives, since $R \in H^1(0,T;L^{\infty}(\Om))$,
\begin{equation} \label{estlp}
\|z\|_{H^1(-T,T;L^2(\Om))} \leq C \|f\|_{L^2(\Omega)}\|R\|_{H^1(0,T;L^{\infty}(\Om))}
\leq C \|f\|_{L^2(\Omega)}.
\end{equation}

In order to apply Theorem~\ref{Carleman} and use the appropriate Carleman estimate, we need a solution of the wave equation that vanishes at time $t=\pm T$. Thus, for $0 < \delta < T$ we take the cut-off function $\theta \in C_0^{\infty}(-T,T)$ such that
\begin{itemize}
\item $0 \leq \theta \leq 1$
\item $\theta(t) = 1$, for all $t \in (-T+\delta, T - \delta)$
\end{itemize}
and we define $~v = \theta z$. Then $v$ satisfies:
\begin{equation}
 \left\{ \begin{array}{rclc}
         v_{tt}  - \textnormal{div}( a \nabla v) + pv  & = & \theta(t)f(x)R_t(x,t) 
         				+ 2\theta_ty_{tt}  + \theta_{tt}y_t  & Q\\
                              v & = & 0 & \Sg \\
            v(0) &=&0& \Om\\
              v_t(0) & = & f(x)R(x,0)  & \Om\\
             v(\pm T) = v_t(\pm T) & = & 0 & \Om.
\end{array}
\right. 
\end{equation}
Take $j \in \{1,2\}$, and let $y$ be the function defined in (\ref{defy}) and $\phi$  the weight function, corresponding  to the point $x_j \in \Om_1$.
Notice that 
\begin{equation}\label{phi1}
 \phi(x,t) \leq \phi(x,0)  \quad  \forall (x,t)\in(0,T)\times\Omega.
\end{equation}

Moreover, by definition of $\rho$ and $\bar a$ (see (\ref{rho}) and the definitions below) we also have 
$$  \frac{|x-x_j|}{\rho(x)} \leq 1 + \frac{|x- y(x)|}{\rho(x)}  \leq 1 + \frac{R_j}{\alpha_j} \leq D_0 $$
and then
\begin{eqnarray}
\phi(x,t)  & =  \bar a\frac{|x-x_j|^2}{\rho^2(x)} - \beta t^2 + M \nonumber\\
	    &  \leq   \bar a D_0^2 - \beta t^2 +M. \label{ineqD}
\end{eqnarray}
Then, by the choice of $~T > D_0 \sqrt{ \frac{a_1}{\beta}}~$ we get 
\begin{equation} \label{inegphiM}
\phi(x,\pm T) < M \leq \phi(x,0). \quad \quad \forall \, x \in \Om
\end{equation}
Thus, taking $\delta$ small enough, it is also true that 
\begin{equation}
\label{phi2}
\phi(x,t) < M \leq \phi(x,0).
\end{equation}
for all $x \in \Om$ and $t\in[-T,-T+\delta] \cup [T-\delta, T]$. 

From now on, $C>0$ will denote a generic constant depending on $\Omega$, $T$, $\beta$, $\theta$, $x_1$, $x_2$, $\delta$, $s_0$ and $\lambda_0$ but independent of $s>s_0$ and $\lambda>\lambda_0$. We will 
occasionally use the notation $\partial_t$ for the time derivative.

As in the  proof of Theorem~\ref{Carleman}, we set $\varphi = e^{\la \phi}$, $w_j = e^{s\varphi}v_j$ and 
$$P_1(w)  =  w_{tt} - a\Delta w + s^2\lambda^2\varphi^2E(\phi)w.$$
 
It is easy to check that 
\begin{equation}
\sum_{j=1}^{2} \left<P_1w_j, \partial_tw_j \right>_{L^2(\Om_j  \times  (0,T))} 
= \frac{1}{2} \int_{\Om} |\partial_t w(0)|^2  + X
\end{equation}
where $X$ is a sum of negligible terms such that
$$
\sum_{j=1}^{2} \left<P_1w_j, \partial_tw_j \right>_{L^2(\Om_j  \times  (0,T))} 
\geq \frac{1}{2} \int_{\Om} |\partial_tw(0)|^2  
        - Cs^2\la^3 \int_0^T \int_\Omega \varphi^3|w|^2.
$$
Since we have $w_t(0) = e^{s\varphi(0)}v_t(0) = e^{s\varphi(0)}f(x)R(x,0)$ and $ |R(x,0)|\geq r$, we get (recall that $Q = \Om \times (-T,T)$, $Q_j = \Om_j \times (-T,T)$ and so on).
$$r^2 \int_{\Om}  e^{2s\varphi(0)} |f|^2  \leq
C \left( \sum_{j=1}^{2} \left<P_1w_j, \partial_tw_j \right>_{L^2(\Om_j \times  (0,T) )} 
+ ~s^2\la^3 \iint_Q \varphi^3|w|^2\right).$$
In order to apply the Carleman estimate (Corollary~\ref{carlem}) we consider both weight functions $\varphi^1$ and $\varphi^2$, corresponding to $x_1$ and $x_2\in \Omega_1$ and  we apply the previous estimates to $w^k_j = e^{s\varphi^k}v_j$ for $j, k = 1,2$ and sum up  the inequalities. We obtain, for $s>s_0$ and $\lambda>\lambda_0$, using Cauchy-Schwarz inequality, the following:
\begin{eqnarray*}
\lefteqn{ r^2 \int_{\Om} ( e^{2s\varphi^1(0)} + e^{2s\varphi^2(0)}  )|f|^2 }\\
&\leq & C \sum_{k=1}^{2}\left( \sum_{j=1}^{2} \left<P^{\phi^k}_1w^k_j, \partial_tw^k_j \right>_{L^2(\Om_j \times (0,T))} 
+ ~s^2\la^3 \iint_Q (\varphi^k)^3|w^k|^2\right)\\
 & \leq & C\sum_{j,k=1}^{2} \left( \frac{1}{\sqrt{s}} \left |P^{\phi^k}_1w^k_j\right|^2_{L^2(Q_j)} 
 +\sqrt{s} \left | \partial_tw^k_j\right|^2_{L^2(Q_j)} \right) \\
   &   &  + ~C~s^2\la^3\sum_{k=1}^2 \iint_Q (\varphi^k)^3\left|w^k\right|^2 \\
& \leq &  \frac{C}{\sqrt{s}}   \sum_{j,k=1}^{2} \left( \left|P^{\phi^k}_1w^k_j\right|^2_{L^2(Q_j)} +
\left\| w^k_j\right\|^2_{Q_j, \varphi^k}    \right) 
\end{eqnarray*}
Now, applying  Corollary~\ref{carlem}, we get
\begin{eqnarray*}
\lefteqn{ r^2 \int_{\Om} ( e^{2s\varphi^1(0)} + e^{2s\varphi^2(0)}  )|f|^2 }\\
& \leq &  \frac{C}{\sqrt{s}}   \sum_{k=1}^{2} \left( \left|P^{\phi^k}_1w^k\right|^2_{L^2(Q_0)} +
\left\| w^k\right\|^2_{Q_0, \varphi^k}    \right)\\ 
& \leq & \frac{C}{\sqrt{s}} \sum_{k =1}^2 \left(
 \left|e^{s\varphi^k}L_pv\right|^2_{L^2(Q_0)} \right.  
    \left.+  s \la  \iint_{\Sigma_+^{\phi^k}} \varphi^k\left |a_2 \dn{w^k}\right |^2  \right)  
\end{eqnarray*}
On the one hand, since we have $\theta_t = 0$  in $[-T+\delta, T-\delta]$, then from estimate (\ref{estlp}), (\ref{phi1}) and (\ref{phi2}), we obtain $\forall k=1,2$
\begin{eqnarray*}
&& \iint_{Q_0} e^{2s\varphi^k}|L_pv|^2 =  \iint_{Q_0} e^{2s\varphi^k} \left|\theta fR_t + 2\theta_ty_{tt}  + \theta_{tt}y_t \right|^2\\
& \leq &C\iint_{Q_0} e^{2s\varphi^k}  |f|^2|R_t |^2
+ C\iint_{Q_0} e^{2s\varphi^k} \left(|\theta_t z_{t} |^2 + |\theta_{tt}z |^2\right)\\
& \leq &C\iint_{Q_0} e^{2s\varphi^k(0)}  |f|^2|R_t |^2
+ C\left(\int_{-T}^{-T+\delta}+\int^{T}_{T-\delta}\right)\int_{\Omega_0} e^{2se^{\la M}} \left(|z_{t} |^2 + |z |^2\right)\\			        
& \leq &  C \|R\|_{H^1(0,T;L^\infty)}\int_{\Om} e^{2s\varphi^k(0)}|f|^2 
+ C e^{2se^{\la M}} \|z\|^2_{H^1(-T,T;L^2(\Om))}\\
& \leq &  C\int_{\Om} e^{2s\varphi^k(0)}|f|^2 
+ C e^{2se^{\la M}} \|f\|_{L^2(\Omega)}\\
& \leq &  C\int_{\Om} e^{2s\varphi^k(0)}|f|^2 .
\end{eqnarray*}
Now, recalling the notation stated in (\ref{bd_0}), we have $\Sigma_+^{\phi^1}\cup\Sigma_+^{\phi^2} \subset \Sigma$ and 
$\left|\dn{w^k}\right| = e^{s\varphi^k}\left|\dn{v}\right| $ on $\Sigma$ for each $k$, so 
we finally obtain
\begin{eqnarray}
\lefteqn{ r^2 \int_{\Om} ( e^{2s\varphi^1(0)} + e^{2s\varphi^2(0)}  )|f|^2 } \nonumber \\
& \leq & \frac{C}{\sqrt{s}} \int_{\Om}\left( e^{2s\varphi^1(0)} + e^{2s\varphi^2(0)} \right )|f|^2 \label{raiz}\\
&   &  + ~C\sqrt{s}\la  \iint_{\Sigma} \left(\varphi^1 e^{2s\varphi^1} + 
              \varphi^2 e^{2s\varphi^2}\right)\left|a_2 \dn{v}\right |^2.  \nonumber
\end{eqnarray}
For $s$ large enough, the left hand side in (\ref{raiz}) can absorb the first term of the right hand side.
Therefore, since $\varphi^k$ and  $\theta$ are bounded on $\Sigma$ and $\dn{z}$ is an even function with respect to $t \in [-T,T]$, we obtain
\begin{eqnarray*}
\int_{\Om} |f|^2 
& \leq & ~C  \iint_{\Sigma} \left(\varphi^1 e^{2s\varphi^1} + 
              \varphi^2 e^{2s\varphi^2}\right)\left|a_2 \dn{v}\right |^2\\
              & \leq & ~C  \iint_{\Sigma} \left|a_2 \dn{z}\right |^2\\
              & = & ~2 C  \int_0^T \int_{\Gamma} \left|a_2 \dn{z}\right |^2.  
\end{eqnarray*}
and this ends the proof of Theorem~\ref{teolp}.
\QED

We will end this paper by the proof of Theorem~\ref{ip} which is a direct consequence of Theorem~\ref{teolp}. Indeed, if  we set $\tilde y = u(q) - u(p)$, $f = p - q$ and $R = u(p)$, then $\tilde y$ is the solution of 
\begin{equation}\label{pbl}
\left\{ \begin{array}{rclc}
\tilde y_{tt}  - \textnormal{div}( a \nabla \tilde y) +( p-f)\tilde y  & = & f(x)R(x,t) & (0,T)\times\Omega\\
\tilde y & = & 0 &(0,T)\times\Sigma \\
\tilde y(0) &=&0& \Om\\
\tilde y_t(0) & = & 0  & \Om
\end{array}
\right. 
\end{equation}
where $q=p-f \in \mathcal U$, with $\mathcal U$ bounded in $L^\infty(\Om)$ from the hypothesis of Theorem~\ref{ip}. The key point is that in the proof of Theorem~\ref{teolp}, all the constants $C>0$ depend on the $L^\infty$-norm of the potential as stated in Corollary \ref{carlem}. Thus, with $q\in \mathcal U$, we are actually, with equation (\ref{pbl}), in a situation similar to the linear inverse problem related to equation (\ref{lp}) and we then obtain the desired result.
\QED

\bigskip\noindent
{\bf Acknowledgments} \\
The authors acknowledge the referees for their useful remarks.

\bigskip

\bigskip\noindent
{\bf References} \\

\end{document}